\documentclass{qam-l}

\def\sq{\hbox {\rlap{$\sqcap$}$\sqcup$}}
\def\C{{\rm C \kern -.15cm \vrule width.5pt \kern .12cm}}
\def\Z{{\rm Z \kern -.27cm \angle \kern .02cm}}
\def\N{ {\rm N \kern -.26cm \vrule width.4pt \kern .10cm}}
\def\1{{\rm 1\mskip-4.5mu l} }
\def\lsim{\raise0.3ex\hbox{$<$\kern-0.75em\raise-1.1ex\hbox{$\sim$}}}
\def\gsim{\raise0.3ex\hbox{$>$\kern-0.75em\raise-1.1ex\hbox{$\sim$}}}

\theoremstyle{definition}

\usepackage{mathrsfs}
\theoremstyle{remark}

\numberwithin{equation}{section}

\begin{document}

\title[Morawetz inequality and asymptotic completeness]{Quadratic Morawetz inequalities and asymptotic completeness in the energy space for nonlinear Schr\"odinger and Hartree equations}

\author{Jean Ginibre}

\address{Laboratoire de Physique Th\'eorique, Universit\'e de Paris XI, B\^atiment
210, F-91405 ORSAY Cedex, France}\thanks {Unit\'e Mixte de Recherche (CNRS) UMR 8627}

\email{Jean.Ginibre@th.u-psud.fr}

\author{Giorgio Velo}

\address{Dipartimento di Fisica, Universit\`a di Bologna and INFN, Sezione di
Bologna, Italy} 
\email{Velo@bo.infn.it}

\dedicatory{\it Dedicated to Professor Walter Strauss on his 70th birthday}


\subjclass[2000]{Primary 35P25; Secondary 35B40,35Q40,35Q55}  

\begin{abstract}
Recently several authors have developed multilinear and in particular quadratic extensions of the classical Morawetz inequality. Those extensions  provide (among other results) an easy proof of asymptotic completeness in the energy space for nonlinear Schr\"odinger equations in arbitrary space dimension and for Hartree equations in space dimension greater than two in the noncritical cases. We give a pedagogical review of the latter results.
\end{abstract}

\maketitle

\section{Introduction}
This paper is devoted to the exposition of an elementary subset of some recent results bearing on scattering theory for the nonlinear Schr\"odinger (NLS) equation
$$i \partial_t u = - (1/2) \Delta u + g(\rho ) u \eqno(1.1)$$

\noindent in $n$ space dimensions, where $u$ is a complex valued function defined in space time ${I\hskip-1truemm R}^{n+1}$, $\rho = |u|^2$ and $g$ is a real valued function of $\rho$, typically a sum of powers 
$$g(\rho ) = \lambda_1 \ \rho^{(p_1-1)/2} + \lambda_2 \  \rho^{(p_2-1)/2} \eqno(1.2)$$

\noindent with $1 < p_1 <p_2$ and $\lambda_1 , \lambda_2 \in {I\hskip-1truemm R}$. We also present a straightforward extension of those results to the Hartree equation (1.1) with 
$$g(\rho ) = V \star \rho \eqno(1.3)$$

\noindent where $V$ is a real even function of the space variable and $\star$ denotes the convolution in ${I\hskip-1truemm R}^n$. The first main question of scattering theory is the existence of the wave operators, namely the construction of solutions that behave asymptotically in time as solutions of the free Schr\"odinger equation, namely such that 
$$u(t) \sim U(t) u_+ \qquad \hbox{for $t \to \infty$} \eqno(1.4)$$

\noindent (and the analogue for $t \to - \infty )$, where 
$$U(t) = \exp (i (t/2)\Delta )\ . \eqno(1.5)$$

\noindent The second main question of scattering theory is asymptotic completeness (AC), simply called ``scattering'' in some of the recent literature, and consists in proving that all solutions of the relevant equation, in a suitable functional framework, behave asympotically as solutions of the free Schr\"odinger equation, namely satisfy (1.4) and the analogue for $t\to - \infty$. Of special interest is the case of finite energy solutions of (1.1), namely of solutions in $L^{\infty}({I\hskip-1truemm R} , H^1)$. An essential tool in the proof of AC for such solutions is the Morawetz inequality, first derived for the nonlinear Klein-Gordon (NLKG) equation \cite{13r} and then extended to the NLS equation \cite{12r}. That inequality was applied to prove AC first for the NLKG equation and then for the NLS equation in space dimension $n \geq 3$ in seminal papers by Morawetz and Strauss \cite{14r} and by Lin and Strauss \cite{12r}, for slightly more regular solutions. The case of general finite energy solutions in space dimension $n \geq 3$ was treated later in \cite{10r} for NLS equations and in \cite{11r} for Hartree equations. The treatment was then improved for the NLS equation in \cite{15r} which covers in addition the more difficult cases $n = 1,2$, as well as the case of the NLKG equation, and for the Hartree equation in \cite{16r}. \par

More recently, several groups of authors have studied the more difficult problem of extending some of the previous results, in particular the proof of existence of global solutions and the proof of AC (``scattering''), on the one hand to the case of critical interactions, and on the other hand to the case of subenergy solutions, namely of solutions of intermediate regularity between $L^2$ and $H^1$.  An important tool in some of those works is a new version of the Morawetz inequality, of a multilinear and in particular bilinear or quadratic type. That inequality has appeared in various forms in the literature and seems to have now stabilized to a simple form (\cite{1r} \cite{3r}-\cite{9r} \cite{14newr} \cite{17r} \cite{19r}-\cite{25r} and references therein quoted). Leaving aside the difficult problems arising for critical interactions and/or for subenergy solutions, that new inequality provides a unified proof of AC for noncritical NLS in the energy space for all space dimensions, as well as for the Hartree equation for $n\geq 3$. That proof is much simpler than the previous ones. The present paper is devoted to an exposition of that new quadratic Morawetz inequality and of its application to the proof of AC for the NLS and Hartree equations in the energy space in noncritical situations. That result for noncritical NLS appears as a by product for $n = 1$ in \cite{4r}, for $n=2$ in \cite{4newr} and for $n \geq 3$ in \cite{20r}, which is mostly devoted to the critical cases. In Section~2, we first derive the quadratic Morawetz identity and we deduce therefrom the basic estimate that leads to the proof of AC. The formal proof of the identity is formulated in terms of conservation laws, in the spirit of \cite{17r}. In Section~3, we exploit the previous estimate to prove AC. We treat the case of the NLS equation in some detail, and we give the modifications needed for the Hartree equation. The result applies to $L^2$ supercritical and $H^1$ subcritical nonlinearities. Some peripheral results are collected in Appendices. In Appendix~1, we give an estimate which points to the usefulness of the Morawetz inequality at lower regularity levels than $H^1$, in particular at the level of $H^{1/2}$. In Appendix~2, we exploit the point of view of conservation laws to derive a quadratic identity for the NLKG equation. That identity however does not lead to estimates because of a lack of positivity. In Appendix~3, we rewrite the original Morawetz inequality for the NLS equation in a form which exhibits its relation to the quadratic identity derived in Section~2. In Appendix~4, we justify the formal computation of Section~2 by a suitable limiting procedure.\par

We conclude this introduction by giving some notation and estimates which will be used freely throughout this paper. For any integer $n \geq 1$, for any $r$, $1 \leq r \leq \infty$, we denote by $\parallel \cdot \parallel_r$ the norm in $L^r \equiv L^r ({I\hskip-1truemm R}^n)$, by $\overline{r}$ the conjugate exponent defined by $1/r + 1/\overline{r} = 1$, and we define $\delta (r) \equiv n/2 - n/r$. We denote by $<\cdot , \cdot >$ the scalar product in $L^2$. We shall use the Sobolev spaces $\dot{H}_r^{\sigma} \equiv \dot{H}_r^{\sigma}({I\hskip-1truemm R}^n)$ and $H_r^{\sigma} \equiv H_r^{\sigma}({I\hskip-1truemm R}^n)$ defined for $0 \leq \sigma < \infty$ and $1 < r < \infty$ by 
$$\dot{H}_r^{\sigma} = \left \{ u : \parallel u, \dot{H}_r^{\sigma}\parallel \ \equiv\ \parallel \omega^{\sigma} u\parallel_r < \infty \right \}$$ 
$$H_r^{\sigma} = \left \{ u : \parallel u, H_r^{\sigma}\parallel \ \equiv\ \parallel <\omega>^{\sigma} u\parallel_r < \infty \right \}$$ 

\noindent where $\omega = (-\Delta )^{1/2}$ and $<\cdot > = (1 + |\cdot |^2)^{1/2}$. The subscript $r$ will be omitted if $r=2$. For any interval $I$ of ${I\hskip-1truemm R}$, for any Banach space $X$, we denote by $\mathscr{C}(I, X)$ the space of continuous functions from $I$ to $X$ and, for $1 \leq q \leq \infty$, by $L^q (I, X)$ (resp. $L_{loc}^q(I,X)$) the space of measurable functions from $I$ to $X$ such that $\parallel u(\cdot ); X \parallel\ \in L^q(I)$ (resp. $\in L_{loc}^q(I)$). \par

We introduce the following definition. A pair of exponent $(q,r)$ is admissible if $0 \leq 2/q = \delta (r) = \delta$ and $\delta \leq 1/2$ for $n=1$, $\delta < 1$ for $n=2$ and $\delta \leq 1$ for $n\geq 3$. Then the well known Strichartz estimates take the form~:\\

\noindent {\bf Lemma 1.1.} {\it Let $U(t)$ be given by (1.5). Then \par

(1) For any admissible pair $(q,r)$
$$ \parallel  U(t) v; L^q ({I\hskip-1truemm R},L^r) \parallel\ \leq \ C  \parallel v \parallel _2 \ .  \eqno(1.6)$$

(2) For any admissible pairs $(q_i,r_i)$, $i=1,2$, and for any interval $I$ of ${I\hskip-1truemm R}$}
$$\parallel \int_{I \cap \{t':t'\leq t \}}  dt' \ U(t-t') \ f(t') ; L^{q_1} (I, L^{r_1})\parallel \ \leq \ C \parallel  f; L^{\overline{q}_2} (I, L^{\overline{r}_2})\parallel  \eqno(1.7)$$

\noindent {\it where the constant $C$ is independent of $I$}. \par \vskip 3 truemm

Lemma~1.1 suggests to study the Cauchy problem for the equation (1.1) in spaces of the following type. Let $I$ be an interval of ${I\hskip-1truemm R}$. We define 
$$X_{(loc)} (I) = \left \{ u : u \in \mathscr{C}(I,L^2)\ {\rm and}\ u \in L_{loc}^q (I, L^r) \hbox{ for all admissible $(q,r)$}\right \}$$

\noindent and
$$X_{(loc)}^1(I) = \left \{ u ; u, \nabla  u \in X_{(loc)} (I) \right \}\ .$$

\section{Quadratic Morawetz inequalities}
In this section we derive the quadratic Morawetz identity for the NLS and Hartree equations and we deduce therefrom the basic estimates that lead to the proof of asymptotic completeness in the energy space for those equations. We begin with a formal derivation of the quadratic Morawetz identity for the NLS equation, assuming sufficient smoothness and decay at infinity of the solutions to give a meaning to the calculation and in particular to the integrations by parts. The underlying algebraic structure is a pair of related conservation laws
$$\partial_t \rho + \nabla \cdot j = 0 \eqno(2.1)$$
$$\partial_t j + \nabla \cdot T = 0 \ . \eqno(2.2)$$

The first one is a scalar conservation law with scalar density $\rho$ and vector current $j$, the second one is a vector conservation law with vector density $j$ and second rank tensor current $T$, and the two laws are related by the fact that the current $j$ of the first one is at the same time the density of the second one. That situation occurs for the NLS equation and with a minor modification for the Hartree equation, as we shall review in this paper. It also occurs for any space time translation invariant system with a symmetric energy momentum tensor, with $\rho$ and $j$ being respectively the energy and momentum densities, and in particular for a class of NLKG equations, for which however it does not lead to useful estimates because of a lack of positivity (see Appendix 2 for the relevant formal calculation in that case). One could also consider the more general situation of two unrelated conservation laws, but that does not seem to be useful in the present case. Let now $h$ be a sufficiently regular real even function defined in ${I\hskip-1truemm R}^{n}$. The starting point is the auxiliary quantity (which will be mostly forgotten at the end)
$$J = (1/2) <\rho , h \star \rho > \ . \eqno(2.3)$$

\noindent From (2.1) (2.2) and with two integrations by parts, it follows that 
$$M \equiv \partial_t J = - \ <\rho, h \star \nabla \cdot j>\ = - \ < \rho , \nabla h \star j > \ , \eqno(2.4)$$
$$\partial_t M = \partial_t^2 J = \ <\nabla \cdot j, \nabla h \star j >\ + \ <\rho, \nabla h \star \nabla \cdot T >$$
$$= - \ <j, \nabla^2 h \star j>\ + \ <\rho , \nabla^2 h \star T> \quad \ \ \eqno(2.5)$$

\noindent where $\nabla^2h$ is the second rank tensor $\nabla_k\nabla_{\ell} h$ and contractions are performed in the obvious way. The quadratic Morawetz identity is then the identity 
$$\partial_t M = - \partial_t \  < \rho , \nabla h \star j>\  = - \  < j, \nabla^2 h \star j>\  + \  <\rho, \nabla^2 h \star T > \ . \eqno(2.6)$$

We now consider the NLS equation
$$i \partial_t u = - (1/2) \Delta u + gu \eqno(1.1)\equiv (2.7)$$

\noindent where $g = g(\rho )$ is a real function of $\rho = |u|^2$. That equation is the Euler-Lagrange equation with Lagrangian density
$$L(u) = \ -\  {\rm Im}\  \overline{u} \partial_t u - (1/2) |\nabla u|^2 - G(\rho ) \eqno(2.8)$$

\noindent where
$$G(\rho ) = \int_0^{\rho} d \rho ' \ g(\rho ')\ . \eqno(2.9)$$

\noindent The basic structure (2.1) (2.2) is realized with $\rho = |u|^2$ and 
$$j = \  {\rm Im} \ \overline{u} \nabla u \ , \eqno(2.10)$$

\noindent and (2.1) is the conservation law of the mass (or charge). The mass current $j$ turns out to be the momentum density, and (2.2) becomes the momentum conservation law. In fact the energy momentum tensor $\widetilde{T}$ is given by 
$$\left \{ \begin{array}{l} \widetilde{T}_{0\ell} = 2{\rm Re} \ \displaystyle{{\partial L \over \partial (\partial_t u)}} \nabla_{\ell} u = - \ {\rm Im} \ \overline{u} \nabla_{\ell} u = - j_{\ell} \\ \\  \widetilde{T}_{k\ell} = 2{\rm Re} \ \displaystyle{{\partial L \over \partial (\nabla_k u)}} \nabla_{\ell} u - \delta_{k\ell} L = -\  {\rm Re} \ \nabla_k \overline{u} \nabla_{\ell} u - \delta_{k\ell}  L \end{array}\right .  \eqno(2.11)$$
\noindent and (2.2) coincides (up to sign) with the conservation law
$$\partial_t \  \widetilde{T}_{0\ell} + \nabla_k \ \widetilde{T}_{k\ell} = 0 \ . \eqno(2.12)$$

\noindent with $T_{k\ell} = - \widetilde{T}_{k\ell}$. For $u$ a solution of (2.7), $L(u)$ reduces to 
$$L(u) = - (1/4) \Delta \rho + \rho g(\rho ) - G(\rho ) \eqno(2.13)$$

\noindent so that
$$T_{k\ell} = \ {\rm Re}\ \nabla_k \overline{u} \nabla_{\ell} u - \delta_{k\ell} \left ( (1/4) \Delta \rho - \rho g + G\right ) \ . \eqno(2.14)$$

\noindent The conservation law (2.2) then holds with $j$ and $T$ defined by (2.10) (2.14), namely
$$\partial_t j = - \nabla \cdot {\rm Re} \nabla \overline{u} \nabla u + \nabla ((1/4) \Delta \rho - \rho g + G) \eqno(2.15)$$

\noindent which can of course be obtained by a direct computation using (2.7). Substituting (2.10) (2.14) or (2.15) into (2.6) yields 
$$\partial_t M =  \  < \rho , \Delta h \star \left ( - (1/4) \Delta \rho + \rho g - G\right ) >\   - \  < j, \nabla^2 h \star j>$$
$$+ \  <\rho, \nabla^2 h \star \nabla  \overline{u} \nabla  u>  \eqno(2.16)$$

\noindent where we have used the symmetry of $\nabla^2h$ to eliminate the real part condition in the last term. On the other hand 
$$<j, \nabla^2 h \star j>\ = \ <\overline{u} \nabla u, \nabla^2 h \star \overline{u} \nabla u>\ - \ <{\rm Re}\ \overline{u} \nabla u, \nabla^2 h \star \ {\rm Re} \ \overline{u} \nabla u>\qquad \qquad$$
$$= \ <\overline{u} \nabla u, \nabla^2 h \star \overline{u} \nabla u>\ - (1/4)\ <\nabla \rho, \nabla^2h \star \nabla \rho > \eqno(2.17)$$

\noindent so that
$$\partial_t M =  (1/2)\  <\nabla \rho , \Delta h \star \nabla \rho > \ + \ <  \rho , \Delta h \star (\rho g - G)>\ + \ R \eqno(2.18)$$

\noindent where we have used the fact that
$$-\ <  \rho , \Delta h \star \Delta \rho >\ = \ <\nabla \rho, \nabla^2h \star \nabla \rho >\ = \ < \nabla \rho, \Delta h \star \nabla \rho >  \eqno(2.19)$$

\noindent by integration by parts, and where
$$R = \ <\overline{u}u, \nabla^2 h \star \nabla  \overline{u} \nabla u >\ - \ < \overline{u} \nabla u , \nabla^2 h \star \overline{u}\nabla u >$$
$$= (1/2) \int dx \ dy \left ( \overline{u}(x) \nabla  \overline{u}(y) - \overline{u}(y)\nabla  \overline{u}(x)\right ) \nabla^2 h (x-y)$$ 
$$ \left ( u(x) \nabla u (y) - u(y) \nabla u (x) \right ) \ . \eqno(2.20)$$

\noindent Integrating (2.18) over time in an interval $[t_1, t_2]$ yields 
$$\int_{t_1}^{t_2} dt \left \{ (1/2) \ < \nabla \rho , \Delta h \star \nabla \rho >\ + \ <   \rho , \Delta h \star (\rho g - G) > + R \right \}$$
$$\left . = - \ < \rho , \nabla h \star\ {\rm Im} \ \overline{u} \nabla u > \right |_{t_1}^{t_2} \ . \eqno(2.21)$$

That identity will yield useful estimates if $\nabla h \in L^{\infty}$ and if $\nabla^2 h$ is nonnegative as a matrix. Under the latter assumption, $R$ is nonnegative, the first term in the integrand is positive, and the second term is nonnegative if $\rho g - G \geq 0$.\par

We now consider a representative situation where the previous formal computations can be made rigorous. We take $h(x) = |x|$, so that 
$$\left \{ \begin{array}{l} \nabla h = |x|^{-1} x \\ \\ \nabla^2 h = |x|^{-1} \left ( \1 - |x|^{-2} x \otimes x \right ) \ ,\ \Delta h = (n-1) |x|^{-1} \quad \hbox{for $n\geq 2$} \\ \\ \nabla^2 h = \Delta h = 2 \delta (x) \quad \hbox{for $n=1$}\ . \end{array}\right . \eqno(2.22)$$
 
\noindent In that case 
$$< \nabla \rho , \Delta h \star \nabla \rho >\ = c\ <\nabla \rho , \omega^{1-n} \nabla \rho >\ = \ c \parallel \rho ; \dot{H}^{(3 - n)/2} \parallel^2 \eqno(2.23)$$ 
$$<  \rho , \Delta h \star (\rho g - G) >\ = c\ <\rho , \omega^{1-n} (\rho g-G) > \eqno(2.24)$$ 

\noindent where $\omega = (- \Delta )^{1/2}$ and $c$ is a constant depending only on $n$ \cite{18r}. \par

We take for $g$ a sum of two powers
$$g(\rho ) = \lambda_1\ \rho^{(p_1-1)/2} + \lambda_2 \  \rho^{(p_2-1)/2} \eqno(1.2) \equiv (2.25)$$

\noindent with $\lambda_1$, $\lambda_2 \in {I\hskip-1truemm R}$, so that
$$G(\rho ) = 2 \lambda_1 \left ( p_1 + 1\right )^{-1} \ \rho^{(p_1+1)/2} +  2 \lambda_2 \left ( p_2 + 1\right )^{-1} \ \rho^{(p_2+1)/2} \eqno(2.26)$$
$$\rho g(\rho ) - G(\rho ) = \lambda_1\ {p_1 - 1 \over p_1 + 1} \ \rho^{(p_1+1)/2} + \lambda_2\ {p_2 - 1 \over p_2 + 1} \ \rho^{(p_2+1)/2} \eqno(2.27)$$

\noindent with $1 \leq p_1 < p_2$. More general $g$ can be easily accomodated. For $H^1$ subcritical powers, the Cauchy problem for NLS is well known to be locally well posed in $X_{loc}^1$ for initial data in $H^1$ and possibly globally well posed \cite{2r}.\\

\noindent {\bf Proposition 2.1.} {\it Let $h(x) = |x|$ and let $g$ be defined by (2.25) with $1 \leq p_1 < p_2$ and $p_2 < 1+4/(n-2)$ for $n \geq 3$. Let $I$ be an interval and let $u \in X_{loc}^1(I)$ be a solution of the NLS equation (2.7). Then \par

(1) The identity (2.21) holds for any $t_1$, $t_2 \in I$. \par

(2) Let in addition $\lambda_1$, $\lambda_2 \geq 0$ (so that $u \in X_{loc}^1({I\hskip-1truemm R}) \cap L^{\infty} ({I\hskip-1truemm R}, H^1)$). Then $u$ satisfies the estimate 
$$\parallel \rho ; L^2({I\hskip-1truemm R}, \dot{H}^{(3-n)/2})\parallel^2 + \int dt\ < \rho , \omega^{1-n} (\rho g - G) > \ \leq\ C \parallel u \parallel_2^3 \ \parallel u ; L^{\infty} ({I\hskip-1truemm R} , H^1)\parallel \eqno(2.28)$$

\noindent In particular $\rho \in L^2({I\hskip-1truemm R} , \dot{H}^{(3-n)/2})$.} \\

\noindent {\bf Sketch of proof.} The proof of Part (1) consists in making the previous formal computation rigorous under the available regularity properties by introducing suitable cut offs and eliminating them by a limiting procedure. This is done in Appendix~4. At this level of regularity, one checks easily that all the terms in the identity are well defined already in the differential form (2.18). Actually by (2.23) and Sobolev inequalities
$$\begin{array}{lll} < \nabla \rho , \Delta h \star \nabla \rho > &\leq \ C\parallel u \parallel_r^2\ \parallel \nabla u \parallel_2^2 &\hbox{with $n/2 - n/r = 1/2$}\\ \\
&\leq \ C \parallel u ; \dot{H}^{1/2}\parallel^2 \ \parallel \nabla u \parallel_2^2 &\hbox{for $n \geq 2$}\ . \end{array} \eqno(2.29)$$

\noindent Similarly, for $g$ a single power $p$
$$\left | < \rho , \Delta h \star (\rho g-G)> \right | \ \leq\ C\parallel u \parallel_{p+1}^{p+1} \ \parallel \omega^{1-n} |u|^2 \parallel_{\infty}$$
$$\leq \ C \parallel u \parallel_{p+1}^{p+1} \ \parallel u \parallel _{r_+} \ \parallel u \parallel_{r_-} \qquad \hbox{for $n \geq 2$} \eqno(2.30)$$

\noindent with $n/r_{\pm} = n/2 - 1/2 \pm \varepsilon$. Furthermore $R = 0$ for $n=1$, while for $n \geq 2$, $R$ is the sum of terms of the type $< \overline{u}\nabla  u, \nabla^2 h \star  \overline{u} \nabla u>$ which are estimated as in (2.29), and $< \overline{u}  u, \nabla^2 h \star  \nabla \overline{u} \nabla u>$ which are estimated by
$$\left | < \overline{u} u, \nabla^2  h \star \nabla \overline{u} \nabla u>\right | \leq C  <|u|^2, \omega^{1-n} | \nabla u|^2>$$
$$\leq \ C\parallel \nabla u \parallel_2^2 \ \parallel \omega^{1-n} |u|^2\parallel_{\infty}\ \leq\ C \parallel \nabla u \parallel_2^2\ \parallel u \parallel_{r_+} \ \parallel u \parallel_{r_-} \ . \eqno(2.31)$$

\noindent Finally the right hand-side of (2.21) is estimated by 
$$\left | < \rho , \nabla h \star {\rm Im}\ \overline{u} \nabla u > \right | \leq \ \parallel \nabla h\parallel_{\infty}\ \parallel u \parallel_2^3\ \parallel \nabla u\parallel_2 \ . \eqno(2.32)$$

\noindent Part (2) follows from (2.21) by taking the limit $t_1 \to - \infty$, $t_2 \to \infty$, from (2.23) (2.24), from the positivity of $\rho g - G$ and of $R$, and from (2.32). \par\nobreak \hfill $\sq$\par

We now sketch briefly some further developments along the previous lines. First the formal computation leading to (2.18) can easily be extended to yield a bilinear Morawetz identity for two solutions of the NLS equation. Actually the identity (2.18) can also be arrived at by applying the original Morawetz identity \cite{12r} to a suitable tensor product of two solutions of (2.7). Let therefore $u_i$, $i = 1,2$, be two solutions of (2.7), let $\rho_i$, $ j_i$, $T_i$ be the associated density, current and tensor $T$, and let $g_i = g(\rho_i)$, $G_i = G(\rho_i)$. We start from
$$J = (1/2)\ <\rho_1 , h \star \rho_2> \eqno(2.33)$$

\noindent so that
$$M \equiv \partial_t J = - (1/2) \left ( <\rho_1 , \nabla h \star j_2> +  <\rho_2 , \nabla h \star j_1>\right ) \ , \eqno(2.34)$$
$$\partial_t  M = \partial_t^2 J = - <j_1, \nabla^2 h \star j_2> + (1/2) \left ( <\rho_1, \nabla^2 h \star T_2> + <\rho_2, \nabla^2 h \star T_1>\right ) \ . \eqno(2.35)$$

\noindent Substituting (2.10) (2.15) into (2.35) and proceeding as before, we obtain 
$$\partial_t  M = (1/2) \Big \{ <\nabla \rho_1 , \Delta h  \star \nabla \rho_2> +  <\rho_1, \Delta h \star (\rho_2 g_2 - G_2)>$$
$$ +  <\rho_2, \Delta  h \star (\rho_1 g_1 - G_1)>\Big \} + R  \eqno(2.36)$$

\noindent where now
$$R = (1/2) \left \{ <\overline{u}_1 u_1, \nabla^2 h \star \nabla \overline{u}_2 \nabla u_2> + (1 \leftrightarrow 2 )\right \} -  <\overline{u}_1 \nabla u_1, \nabla^2 h \star \overline{u}_2 \nabla u_2>$$
$$= (1/2) \int dx\ dy \left ( \overline{u}_1(x) \nabla \overline{u}_2(y) - \overline{u}_2 (y) \nabla \overline{u}_1(x) \right ) \nabla^2 h (x-y)$$
$$\times \left ( u_1(x) \nabla u_2(y) - u_2 (y) \nabla u_1(x)\right ) \ .\eqno(2.37)$$

\noindent The identity (2.36) is the bilinear version of (2.18). Remarkably enough, $R$ is still nonnegative in that case if $\nabla^2h$ is a nonnegative matrix. On the other hand for a repulsive (defocusing) $g$, the terms in (2.36) containing $g$ are also nonnegative, while the first term in the bracket is positive for $n\geq 3$ (but in general not for $n = 1,2$), so that in that case (2.36) yields some bilinear estimates. Whether such estimates can be useful remains to be seen.\par

A second further development of the previous calculation consists in using for $h$ other functions than $|x|$. For instance one can take $h(x) = |\theta \cdot x|$ for $\theta \in S^{n-1}$ and more generally $h(x) = |Px|$ for $P$ the orthogonal projection on a generic $k$ dimensional plane in ${I\hskip-1truemm R}^n$. The first choice leads naturally to an estimate of the Radon transform of $\rho$ \cite{17r}. One can also take advantage of the fact that the derivation of (2.18) involves mainly two integrations by parts from $h$ to $\nabla^2h$ in order to treat the case of a domain $\Omega \subset {I\hskip-1truemm R}^n$, typically the complement of a convex (or at least star-shaped) compact subset of ${I\hskip-1truemm R}^n$. One then obtains identities similar to (2.18) with additional surface terms, from which one can derive estimates of solutions in $\Omega$ \cite{17r}. \par

A third possible development consists in extending the estimates of Proposition 2.1, part (2) to the case of attractive (focusing) interactions $g$ and of small solutions. We consider for illustration the case of a single power
$$g(\rho ) = - \rho^{(p-1)/2} \eqno(2.38)$$

\noindent in dimension $n=1$. In that case, (2.18) becomes (remember that $R = 0$ for $n=1$)
$$\partial_t  M = \ \parallel  \nabla \rho \parallel_2^2 - {p-1 \over p+1} \int dx \ \rho^{(p+3)/2}\ . \eqno(2.39)$$

\noindent By Sobolev inequalities, we estimate 
$$\parallel  \rho \parallel_{(p+3)/2}^{(p+3)/2} \ \leq \ C  \parallel  \nabla \rho \parallel_2^2 \ \parallel  \rho \parallel_{(p-1)/4}^{(p-1)/2}\ \leq \ C  \parallel  \nabla \rho \parallel_2^2 \ \parallel u; \dot{H}^{\sigma_c}  \parallel^{p-1}\eqno(2.40)$$

\noindent where $\sigma_c = 1/2 - 2/(p-1)$ is the value of $\sigma$ for which $g$ given by (2.38) is $\dot{H}^{\sigma}$ critical, provided $\sigma_c \geq 0$, namely $p \geq 5$. Therefore
$$\partial_t  M \geq\   \parallel  \nabla \rho \parallel_2^2 \left ( 1 - C \parallel u ; \dot{H}^{\sigma_c} \parallel^{p-1}\right )  \eqno(2.41)$$

\noindent so that (2.18) again yields an a priori estimate of $\rho$ in $L^2( {I\hskip-1truemm R}, \dot{H}^1)$ provided $M$ is controlled and provided $u$ is small in $L^{\infty}({I\hskip-1truemm R}, \dot{H}^{\sigma_c})$. The latter condition can be realized for energy solutions by taking some initial data $u_0$ small in $L^2$ if $\sigma_c = 0$, namely $p = 5$, and $u_0$ small in $H^1$ if $p > 5$. That smallness condition is of the same type as that occurring in the proof of boundedness of the $H^1$ norm from the energy conservation law which is used in the standard proof of globalization in $H^1$. \par

A similar situation can occur in higher space dimensions in so far as one can prove the estimate
$$\left | <\rho , \omega^{1-n} \rho^{(p+1)/2}>\right | \leq \ C \parallel \rho ;  \dot{H}^{(3-n)/2}  \parallel^2 \  \parallel u; \dot{H}^{\sigma_c}  \parallel^{p-1}\eqno(2.42)$$

\noindent where again $\sigma_c = n/2 - 2/(p-1)$ is the critical Sobolev exponent corresponding to $p$, provided $\sigma_c \geq 0$, namely $p \geq 1 + 4/n$, the $L^2$ critical value. The estimate (2.42) can be proved easily by the use of Sobolev inequalities for $n = 2,3$ and $p$ not too large. We leave the investigation of that estimate for general $n$ and $p$ as an open question.\par

A last possible development consists in using the Morawetz inequality to prove global wellposedeness and possibly AC (``scattering'') at a lower level of regularity than $H^1$, and that possibility has been extensively exploited. See for instance \cite{3r}-\cite{4r} \cite{6r}-\cite{9r} \cite{16newr} \cite{25r} and references therein quoted. In particular the right hand-side of (2.21) is controlled by the $H^{1/2}$ norm of $u$. For completeness we give a proof of that fact in Appendix 1 (see also \cite{6r} for the case $n \geq 3$). \par

We now turn to the Hartree equation (1.1) with $g$ given by (1.3). The formal computation is almost the same as for the NLS equation, except for the fact that, because of the nonlocality of the interaction, the equation is not Lagrangian. However the evolution equation of $j$ for the NLS equation takes the form 
$$\partial_t  j = \ {\rm kinetic \ terms}\ - \rho \nabla g \eqno(2.43)$$

\noindent as follows in the same way as (2.15) from a computation which can be done without referring to the special form of $g$, so that (2.43) also holds for the Hartree equation (1.1) (1.3). Substituting (2.43) into $\partial_t  M$ and using the fact that the kinetic terms are unchanged, we obtain 
$$\partial_t  M = (1/2) \ <\nabla \rho , \Delta h \star \nabla \rho > + <  \rho , \nabla h \star (\rho \nabla (V \star \rho ))> + \ R \eqno(2.44)$$

\noindent where $R$ is given by (2.20) as before. Integrating (2.44) over time in an interval $[t_1, t_2]$ yields
$$\int_{t_1}^{t_2} dt \left \{ (1/2) \ <\nabla \rho , \Delta h \star \nabla \rho >  +  < \rho, \nabla h \star (\rho \nabla (V \star \rho ))> + \ R \right \} $$
$$\left . = -\  <\rho , \nabla h \star \ {\rm Im} \ \overline{u} \nabla u >\right |_{t_1}^{t_2} \ . \eqno(2.45)$$

\noindent As in the case of the NLS equation, that identity will yield useful estimates if $\nabla h \in L^{\infty}$ and if $\nabla^2 h$ is nonnegative as a matrix, so that $R$ is nonnegative, and if in addition the potential term in (2.45) is nonnegative. We now show that this is the case if $V$ is radial and nonincreasing. Assuming sufficient smoothness and decay at infinity for $V$, we obtain 
$$P \equiv \ <\rho , \nabla h \star (\rho \nabla (V  \star \rho ))> \ = \int dx\ dy\ dz \ \rho (x) \nabla h (x-y) \rho (y) \nabla V (y-z) \rho (z)$$
$$= (1/2) \int dx\ dy\ dz\ \rho (x) \ \rho (y) \ \rho (z) \ \nabla V(y-z) (\nabla h(x-y) - \nabla h (x-z)) \eqno(2.46)$$

\noindent where we have used the fact that $\nabla V$ is an odd function. In order to prove the positivity of that integral, it suffices to prove that for all $x$, $y$
$$\nabla V(x) \cdot (\nabla h (x+y) - \nabla h (y)) \leq 0 \eqno(2.47)$$

\noindent where we have changed variables from $(y-z, x-y, x-z)$ to $(x,y,x+y)$. Let now $V(x) = v(|x|)$. The left hand-side of (2.47) can be written as
$$\int_0^1 d\theta \ \nabla V (x) x \cdot \nabla^2 h (y + \theta x)$$
$$= \int_0^1 d\theta \ |x|^{-1} v'(|x|) (x \otimes x) \cdot \nabla^2 h (y + \theta x) \leq 0 \eqno(2.48)$$

\noindent for nonpositive $v'$ and nonnegative $\nabla^2 h$. \par

We now give a proposition where we assume sufficient regularity of $V$ to ensure wellposedness in $H^1$ and to make the previous formal computation rigorous.\\

\noindent {\bf Proposition 2.2.} {\it Let $h = |x|$ and let $V \in L^{p_1} + L^{p_2}$ where
$$p_2 \geq 1 \quad , \qquad n/4 < p_2 <p_1 \leq \infty\ . \eqno(2.49)$$

\noindent Let $I$ be an interval and let $u \in X_{loc}^1(I)$ be a solution of the Hartree equation (1.1) (1.3). Then\par

(1) The identity (2.45) holds for any $t_1, t_2 \in I$\par

(2) Let in addition $V$ be radial non increasing (so that $V$ is non negative, possibly up to a harmless constant, and $u \in X_{loc}^1({I\hskip-1truemm R}) \cap L^{\infty}({I\hskip-1truemm R}, H^1)$). Then $u$ satisfies the estimate}
$$\parallel \rho ; L^2({I\hskip-1truemm R}, \dot{H}^{(3-n)/2}) \parallel^2\ \leq \ C \parallel u \parallel_2^3 \ \parallel u; L^{\infty}({I\hskip-1truemm R}, H^1) \parallel \ . \eqno(2.50)$$
\vskip 3 truemm

\noindent {\bf Sketch of proof.} The proof of Part (1) follows the same pattern as that of Proposition~2.1. Here we simply verify that the Hartree potential term $P$ in (2.45) is well defined at the available level of regularity. By the H\"older and Young inequalities, we estimate 
$$|P| \equiv \left | <\rho , \nabla h \star (\rho (V \star \nabla \rho ))> \right | \leq \ \parallel \rho \parallel_1 \ \parallel \nabla h \parallel_{\infty}\ \parallel \rho \parallel_{k/2}\ \parallel V \parallel_p \ \parallel u \parallel_k \ \parallel \nabla u \parallel_2$$
$$= \  C \parallel \rho \parallel_1\ \parallel \nabla u \parallel_2 \ \parallel u \parallel_k^3 \eqno(2.51)$$

\noindent with $\delta (k) \equiv n/2 - n/k = n/(3p)$. For the relevant values of $p$, one can take $\delta (k) \leq 1/4$ for $n=1$, $\delta (k) < 1$ for $n = 2$, $\delta (k) \leq 1$ for $n=3$ and for $n \geq 4$ if $p \geq n/3$, so that $ \parallel u  \parallel_k$ is controlled by the $H^1$ norm of $u$ and $P$ is controlled in $L_{loc}^{\infty} (I)$, namely at the differential level. For $n \geq 4$ and $n/4 \leq p < n/3$, we use the fact that $u \in L_{loc}^q(I,L^k)$ with $2/q = \delta (k) - 1 = n/(3p)-1 \leq 1/3$, so that $u \in L_{loc}^6(I,L^k)$ and therefore $P \in L_{loc}^2 (I)$.\par

Part (2) follows from (2.45) by taking the limit $t_1 \to - \infty$, $t_2 \to \infty$, from (2.23) (2.32) and from the positivity of $P$ defined in (2.46) and of $R$ defined by (2.20). \par\nobreak \hfill $\sq$\par

\section{Asymptotic completeness in the energy space}
In this section we exploit the Morawetz estimates of Propositions 2.1 and 2.2 to derive asymptotic completeness in $H^1$ for the NLS and Hartree equations. We begin with the NLS equation for which we restrict our attention to a single power interaction
$$g(\rho ) = \lambda \ \rho^{(p-1)/2} \eqno(3.1)$$

\noindent We shall use the parameter $\sigma_c$ defined equivalently by
$$\sigma_c = n/2 - 2/(p-1) \qquad {\rm or}\quad p-1 = 4/(n-2 \sigma_c) \eqno(3.2)$$

\noindent so that $g$ given by (3.1) is $\dot{H}^{\sigma_c}$ critical. We shall assume $0 < \sigma_c < 1$ so that $g$ is $L^2$ supercritical and $H^1$ subcritical. The treatment extends in a trivial way to a sum of such powers and to more general $g$. The case of critical powers is much more complicated and we refer to \cite{20r} for a treatment of that case in dimension $n \geq 3$. Some of the arguments can also be applied to solutions in $H^{\sigma}$ for $0 < \sigma \leq 1$.\par

The main tehnical step is the following proposition.\\

\noindent {\bf Proposition 3.1.} {\it Let $g$ be defined by (3.1) with $0 < \sigma_c < 1$ ($\sigma_c < 1/2$ for $n = 1$). Let $u \in X_{loc}^1 ({I\hskip-1truemm R}) \cap L^{\infty} ({I\hskip-1truemm R}, H^1)$ be a solution of the NLS equation (1.1) (3.1) such that $\rho = |u|^2 \in L^2 ( {I\hskip-1truemm R}, \dot{H}^{(3-n)/2})$. Then $u \in X^1 ({I\hskip-1truemm R})$.}\\

\noindent {\bf Remark 3.1.} For repulsive (defocusing) interaction $g$, namely for $\lambda > 0$, the Cauchy problem with initial data in $H^1$ is known to yield solutions satisfying the first assumption, and those solutions satisfy the condition on $\rho$ by Proposition 2.1. For attractive (focusing) interaction, the first assumption is satisfied for small data in $H^1$, and the assumption on $\rho$ can also be satisfied in some cases, for instance for $n = 1$, and for $n = 2,3$ and $p$ not too large, as discussed in the comments after Proposition 2.1 (see in particular (2.41)).\\

\noindent {\bf Proof.} Let $I = [t_0, t_1]$ be an interval and $u_0 = u(t_0)$. We start from the integral equation
$$u(t) = U(t-t_0)u_0 - i \int_{t_0}^t dt' \ U(t-t') \ g(\rho (t')) \ u(t')\ . \eqno(3.3)$$

\noindent  Using the Strichartz inequalities, we estimate in a standard way \cite{2r}
$$\parallel u; X^1(I)\parallel \ \leq \ C\left ( \parallel u_0:H^1\parallel\ + \ \parallel g(\rho ) u; L^{\overline{q}} (I, H_{\overline{r}}^1)\parallel \right )$$
$$\leq C \left ( \parallel u_0;H^1\parallel\ + \ \parallel u; X^1(I)\parallel\ \parallel u; L^k(I,L^{\ell})\parallel^{p-1} \right ) \eqno(3.4)$$

\noindent where $1/\overline{r} + 1/r = 1/\overline{q} + 1/q = 1$, $(q,r)$ is an admissible pair, and 
$$\left \{ \begin{array}{l}2/k = \left ( n/2 - \sigma_c\right  ) (1 - \delta )\\ \\ n/\ell = \left ( n/2 - \sigma_c \right ) \delta \end{array} \right . \eqno(3.5)$$

\noindent where $\delta \equiv \delta (r) = n/2 - n/r$. \par

The main step of the proof consists in estimating $u$ in $L^k(L^{\ell})$ by interpolation between the Morawetz quantity $\parallel \rho ; L^2(\dot{H}^{(3-n)/2}) \parallel$ and some norm which is controlled by\break\noindent $\parallel u; L^{\infty} (H^1)\parallel$, typically $\parallel u; L^{\infty} (\dot{H}^{\sigma})\parallel$ for some $\sigma$, $0 \leq\sigma \leq 1$. For orientation, we first consider the homogeneity degree of the various norms involved, where the degree of $\parallel u ; L^q (\dot{H}^{\sigma}_r)\parallel$ is defined as $\sigma + \delta (r) - 2/q$, so that it reduces to $\sigma$ for admissible $(q,r )$. In particular the degree of $\parallel u ;L^k(L^{\ell})\parallel$ is $\sigma_c$ by (3.5), that of $\parallel u ; L^{\infty} (\dot{H}^{\sigma})\parallel$ is $\sigma$ and the degree $\sigma_M$ of the Morawetz quantity is obtained by comparing from the point of view of dimension 
$$\parallel \rho ; L^2 (\dot{H}^{(3-n)/2})\parallel \ \sim  \ \parallel u ; L^{\infty} ( \dot{H}^{\sigma_M})\parallel^2\ , $$

\noindent which gives 
$$1 + n/2 + (n-3)/2 = 2\left ( n/2 - \sigma_M\right )$$

\noindent and therefore $\sigma_M = 1/4$. \par

We have to combine information on $u$ and on $\rho$, which can be transformed into information bearing only on $u$ or only on $\rho$. We consider separately the cases $n \geq 2$ and $n = 1$. \\

\noindent {\bf The case n $\geq$ 2.} Here we work with $u$. The information on $\rho$ implies the following information on $u$.

$$\left \{ \begin{array}{ll}\rho \in L^2 ({I\hskip-1truemm R}, \dot{H}^{1/2}) \subset L^2 ({I\hskip-1truemm R}, L^4) \Leftrightarrow u \in L^4 ({I\hskip-1truemm R}, L^8)  &\hbox{for $n = 2$} \\ \\  \rho \in L^2 ({I\hskip-1truemm R}, L^{2})\Leftrightarrow u \in L^4 ({I\hskip-1truemm R}, L^4) &\hbox{for $n = 3$} \\ \\  \rho \in L^2 ({I\hskip-1truemm R}, \dot{H}^{(3-n)/2})  \Rightarrow u \in L^4 ({I\hskip-1truemm R}, \dot{H}_4^{(3-n)/4}) &\hbox{for $n \geq 4$} \end{array}\right . \eqno(3.6)$$

\noindent where the last result follows from Lemma 5.6 in \cite{20r}. We want to estimate $u$ in $L^k (I, L^{\ell})$ with $k$, $\ell$ satisfying (3.5) for some $k < \infty$ and some $\delta$ with $0 \leq \delta < 1$ (the value $\delta = 1$ is excluded a priori for $n=2$, and by the condition $k < \infty$ for $n \geq 3$). From (3.5), we obtain 
$$2/k + n/\ell = n/2 - \sigma_c \ . \eqno(3.7)$$

\noindent Conversely if $k$, $\ell$ satisfy (3.7) with $1 \leq k < \infty$ and $2 \leq \ell \leq \infty$, then $\delta$ defined by (3.5) satisfies $0 \leq \delta < 1$, so that it suffices to consider (3.7). We estimate by Sobolev inequalities and by (3.6)
$$\parallel u; L^k(I, L^{\ell} )\parallel \ \leq \ C\ \parallel u;L^4(I, \dot{H}_4^{(3-n)/4}) \parallel^{\theta}  \ \parallel u; L^{\infty} (I, \dot{H}^{\sigma })\parallel^{1 - \theta} \eqno(3.8)$$

\noindent (where $\dot{H}_4^{1/4}$ should be replaced by $L^8$ for $n = 2$ according to (3.6)) for some $\sigma$ and $\theta$ with $0 \leq \sigma \leq 1$ and $0 < \theta \leq 1$, such that 
$$\left \{ \begin{array}{l} 2/k = \theta /2 \\ \\ n/\ell = \theta (n/2 - 3/4) + (1 - \theta ) (n/2 - \sigma )  \end{array}\right . \eqno(3.9)$$

\noindent so that
$$2/k + n/\ell = n/2 - \sigma_c = \theta (n/2 - 1/4) + (1 - \theta ) (n/2 - \sigma )$$

\noindent or equivalently
$$\sigma_c = \theta /4 + (1 - \theta ) \sigma \eqno(3.10)$$

\noindent in accordance with the homogeneity argument given above. In addition for $n \geq 4$, the Sobolev inequality requires
$$\theta (n-3)/4 \leq (1 - \theta ) \sigma \ . \eqno(3.11)$$

\noindent For a given $\sigma_c$ with $0 < \sigma_c < 1$, it is therefore sufficient to find $\sigma$ and $\theta$ with $0 \leq \sigma \leq 1$ and $0 < \theta \leq 1$, satisfying (3.10) and in addition (3.11) or equivalently
$$\theta \leq 4\sigma /(n-3 + 4 \sigma )   \eqno(3.12)$$

\noindent for $n \geq 4$. One can make the following choices.\\

\noindent {\bf Case n = 2, 3.} For $\sigma_c = 1/4$, one can take $\theta = 1$ and the norm in $L^{\infty} (\dot{H}^{\sigma})$ is not needed. For $\sigma_c \not= 1/4$, the allowed values of $\sigma$ are given by
$$0 \leq \sigma < \sigma_c < 1/4\qquad \hbox{or} \quad 1/4 < \sigma_c  < \sigma \leq 1\ , \eqno(3.13)$$

\noindent with $\theta$ defined by (3.10).\\

\noindent {\bf Case n $\geq$ 4.}  For $\sigma_c = 1/4$, one must  take $\sigma = 1/4$ and one can take $\theta = (n-2)^{-1}$. For $\sigma_c \not= 1/4$, the allowed values of $\sigma$ are given by 
$$(0 <) \sigma_0 \leq \sigma < \sigma_c  < 1/4\qquad \hbox{or} \quad 1/4 < \sigma_c  < \sigma \leq \sigma_0 \wedge 1\ , \eqno(3.14)$$

\noindent where $\sigma_0$ is defined by (3.10) and (3.12) with equality, namely
$$\sigma_0 = \sigma_c (n-3) /\left ( n-2-4 \sigma_c\right ) \ , \eqno(3.15)$$

\noindent with $\theta$ defined by (3.10).\par

We can now complete the proof of the proposition. Substituting (3.8) into (3.4) yields
$$\parallel u; X^1(I)\parallel \left ( 1 - M_1 \parallel \rho ; L^2(I, \dot{H}^{(3-n)/2}) \parallel^{\theta (p-1)/2} \right ) \leq M_2 \eqno(3.16)$$

\noindent where $M_1$, $M_2$ depend only on $\parallel u ; L^{\infty} ({I\hskip-1truemm R}, H^1)\parallel$. By Proposition 2.1, one can partition ${I\hskip-1truemm R}$ into a finite number of intervals such that 
$$M_1 \parallel \rho ; L^2(I, \dot{H}^{(3-n)/2}) \parallel^{\theta (p-1)/2}\  \leq 1/2 \eqno(3.17)$$

\noindent and the number of intervals is also estimated in terms of $\parallel u; L^{\infty} ({I\hskip-1truemm R}, H^1)\parallel$. This yields an estimate of $\parallel u; X^1(I)\parallel$ for each such interval. Furthermore $u \in X^1({I\hskip-1truemm R})$ and $\parallel u; X^1({I\hskip-1truemm R})\parallel$ is estimated by a (computable) power of $\parallel u; L^{\infty} ({I\hskip-1truemm R}, H^1)\parallel$. This completes the proof for $n \geq 2$.\\

\noindent {\bf Remark 3.2.} For $n=2,3$, the argument is the same whether one uses $u$ or $\rho$. For $n \geq 4$, the argument can also be made by using $\rho$ and the fact that
$$\parallel \rho ; \dot{H}_{n/(n- \sigma )}^{\sigma} \parallel \ \leq \ C \parallel  u; \dot{H}^{\sigma}\parallel ^2 \eqno(3.18)$$

\noindent for $0 \leq \sigma \leq 1$ by Leibniz and Sobolev inequalities, and one ends up again with the condition (3.14) with however 
$$\sigma_0 = 2 \sigma_c (n-3)/\left ( 2n-5-4 \sigma_c \right ) \eqno(3.19)$$

\noindent which makes the restriction on $\sigma$ slightly stronger. \\

\noindent {\bf The case n = 1.} Here we work with $\rho$. For low values of $p$, we shall need the implication for $\rho$ of some Strichartz norms of $u$. We need the following lemma. \\

\noindent {\bf Lemma 3.1.} {\it Let $0 \leq \sigma < 1/r \leq 1/2$. Then
$$\parallel \rho ; \dot{H}_{(2/r- \sigma )^{-1}}^{\sigma} \parallel \ \leq \ C \parallel  u; \dot{H}_r^{\sigma}\parallel ^2 \eqno(3.20)$$

\noindent and therefore for $2/q = \delta (r)$ and for any interval $I$}
$$\parallel \rho ; L^{q/2} (I, \dot{H}_{(2/r- \sigma )^{-1}}^{\sigma}) \parallel \ \leq \ C \parallel  u; L^q (I, \dot{H}_r^{\sigma})\parallel ^2 \ . \eqno(3.21)$$
\vskip 5 truemm

\noindent {\bf Proof of Lemma 3.1.} We estimate by fractional Leibniz and Sobolev inequalities 
$$\parallel \omega^{\sigma} \rho \parallel_{(2/r - \sigma )^{-1}} \ \leq \ C \parallel \omega^{\sigma} u \parallel_r \  \parallel u  \parallel_{(1/r - \sigma )^{-1}} \ \leq \ C  \parallel \omega^{\sigma} u  \parallel_r^2\ . $$
\hfill $\sq$\par

We come back to the proof of the proposition. We start again from (3.4), so that we need to estimate $u$ in $L^k (I, L^{\ell})$ with 
$$\left \{ \begin{array}{l} 2/k = (1/2 - \sigma_c ) (1 - \delta )\\ \\ 1/\ell =(1/2 - \sigma_c ) \delta  \end{array}\right . \eqno(3.22)$$

\noindent for some $\delta$ with $0 \leq \delta \leq 1/2$, or equivalently with
$$2/k + 1/\ell = 1/2 - \sigma_c \ , \eqno(3.23)$$
$$0 \leq 1/\ell \leq (1/2 - \sigma_c )/2 \ . \eqno(3.24)$$

\noindent We estimate 
$$\parallel u; L^k(I, L^{\ell} )\parallel^2 \ = \ \parallel \rho ; L^{k/2} (I, L^{\ell /2}) \parallel \ \leq \ C\ \parallel \rho ;L^2(I, \dot{H}^1) \parallel^{\theta}  \ \parallel \rho; L^{q/2} (I, \dot{H}_{(2/r - \sigma )^{-1}}^{\sigma })\parallel^{1 - \theta}  \eqno(3.25)$$

\noindent by Sobolev inequalities, for some $\sigma$, $\theta$ and admissible $(q,r)$ satisfying $0 \leq \sigma < 1/r \leq 1/2$, $0 < \theta \leq 1$ and
$$\left \{ \begin{array}{l} 2/k = \theta /2 + (1 - \theta )2/q \\ \\ 1/\ell =-  \theta /4+ (1 - \theta ) (1/r - \sigma ) \ .  \end{array}\right .  \eqno(3.26)$$

\noindent Substituting (3.26) into (3.23) (3.24) yields
$$\sigma_c = \theta /4 + (1 - \theta )\sigma \ ,  \eqno(3.27)$$
$$0 \leq - \theta /4 + (1 - \theta) (1/r-  \sigma ) \leq (1/2 -  \sigma_c )/2 \ .  \eqno(3.28)$$

\noindent For $\sigma_c = 1/4$, namely $p = 9$, we must take $ \sigma = 1/4$ and we ensure (3.27) (3.28) by taking $r = 2$ and $\theta = 1/2$. \par

\noindent For $\sigma_c \not= 1/4$, we must take
$$0 \leq \sigma < \sigma_c < 1/4 \qquad {\rm or}\qquad 1/4 <  \sigma_c <  \sigma < 1/2 \eqno(3.29)$$

\noindent and the elimination of $\theta$ betwen (3.27) (3.28) yields
$$ \sigma_c \leq {(4 \sigma_c - 1) \over (4  \sigma - 1)r} \leq (1 + 2 \sigma_c)/4  \eqno(3.30)$$

\noindent which implies the condition $ \sigma < 1/r$ since
$${1 \over r} \geq {4 \sigma - 1 \over 4  \sigma_c -1} \  \sigma_c =  \sigma + { \sigma -  \sigma_c \over 4  \sigma_c - 1} >  \sigma \eqno(3.31)$$

\noindent by (3.29). One can fulfill (3.30) with $r=2$ provided 
$$ \sigma_+   \ \mathrel{\mathop >_{<}}\  \sigma  \  \mathrel{\mathop >_{<}}\ \sigma_-   (\mathrel{\mathop >_{<}}\   \sigma_c) \ \mathrel{\mathop >_{<}}\ 1/4 \eqno(3.32)$$

 \noindent where
$$ \sigma_+ = \left ( 6 \sigma_c - 1\right )/8 \sigma_c\ , \qquad  \sigma_- = (10  \sigma_c - 1) /(8  \sigma_c + 4)\ , \eqno(3.33)$$

\noindent which is compatible with (3.29) provided $ \sigma_- \geq 0$, namely $ \sigma_c \geq 1/10$ or $p \geq 6$. In that case we take $r=2$ and we can take 
$$\left \{ \begin{array}{ll} 0 \leq \sigma \leq \sigma_- (< \sigma_c) &\hbox{for $1/10 \leq \sigma_c \leq 1/6$} \\ \\ (0 \leq ) \sigma_+ \leq \sigma \leq \sigma_- (< \sigma_c) &\hbox{for $1/6 \leq \sigma_c < 1/4$} \\ \\   \sigma_+ = \sigma = \sigma_- = \sigma_c= 1/4 &\hbox{for $\sigma_c = 1/4$} \\ \\ (\sigma_c <)  \sigma_- \leq \sigma \leq \sigma_+ (< 1/2) &\hbox{for $1/4 < \sigma_c < 1/2$} \end{array}\right .  \eqno(3.34)$$

\noindent with $\theta$ defined by (3.27) for $\sigma_c \not= 1/4$. For such $(\sigma , \theta )$, one obtains 
$$\parallel u ; L^k(I, L^r)\parallel\ \leq \ C\parallel \rho ; L^2(I, \dot{H}^1)\parallel^{\theta /2} \ \parallel u ; L^{\infty}(I, \dot{H}^{\sigma})\parallel^{1 - \theta} \eqno(3.35)$$

\noindent which implies (3.16) with $n=1$. \par

For $0 < \sigma_c < 1/10$, namely $5 < p < 6$, one can take $\sigma = 0$ and take for $r$ the minimal value allowed by (3.30), namely
$$4/r = (1 + 2\sigma_c )/(1 - 4\sigma_c ) \eqno(3.36)$$ 

\noindent and $\theta = 4 \sigma_c$. One then obtains
$$\parallel u ; L^k(I, L^\ell)\parallel\ \leq \ C\parallel \rho ; L^2(I, \dot{H}^1)\parallel^{\theta /2} \ \parallel u ; X(I) \parallel^{1 - \theta} \eqno(3.37)$$

\noindent so that by (3.4)
$$\parallel u ; X^1(I)\parallel\ \leq \ C\Big (  \parallel u ; L^{\infty} (I, H^1)  \parallel \ + \ \parallel \rho ; L^2(I, \dot{H}^1)\parallel^{(p-1) \theta /2}$$
$$\times \ \parallel u ; X^1(I)\parallel^{1 + (p - 1)(1  - \theta)}\Big ) \eqno(3.38)$$

\noindent  which gives again an estimate of $\parallel u, X^1(I) \parallel$ provided $\parallel \rho ; L^2(I, \dot{H}^1)\parallel$ is sufficiently small. \par

The end of the proof proceeds as in the case $n \geq 2$. \par \hfill $\sq$ \par

\noindent {\bf Remark 3.3.} If one wants to use values of $\sigma$ arbitrarily close to $\sigma_c$ for $\sigma_c > 1/4$, one needs to take $r > 2$ in the region $\sigma_c < \sigma < \sigma_-$. The lowest possible value of $r$ is given by (see (3.30))
$$4/r = (1 + 2\sigma_c) (4\sigma - 1) /(4\sigma_c - 1) \ . \eqno(3.39)$$
\vskip 5 truemm

\noindent {\bf Remark 3.4.} One could use $u$ instead of $\rho$ also in the case $n=1$. From the inequality
$$\rho^{3/2} \leq (3/4) \int dx \ \rho^{1/2} |\rho '| \leq (3/4) \parallel \rho \parallel_1^{1/2}\ \parallel \rho '\parallel_2$$

\noindent we obtain
$$\parallel u; L^6 ({I\hskip-1truemm R}, L^{\infty})\parallel^3\ = \ \parallel \rho ; L^3({I\hskip-1truemm R} , L^{\infty})\parallel^{3/2} \ \leq (3/4) \parallel u \parallel_2\ \parallel \rho ; L^2({I\hskip-1truemm R} , \dot{H}^1)\parallel \eqno(3.40)$$

\noindent and one can perform the estimates by using $\parallel u; L^6 ({I\hskip-1truemm R}, L^{\infty})\parallel$ instead of $\parallel \rho : L^2({I\hskip-1truemm R} , \dot{H}^1)\parallel$. The results are essentially the same with however stronger restrictions on $\sigma$. \\

We now exploit Proposition 3.1 to prove AC in $H^1$ for the NLS equation (1.1) with interaction (3.1). We first recall some standard results on scattering for that equation \cite{2r}.\\

\noindent {\bf Proposition 3.2.} {\it Let $0 \leq \sigma_c < 1$, $\sigma_c < 1/2$ for $n=1$, or equivalently $p \geq 1 + 4/n$, $p < 1 + 4/(n-2)$ for $n \geq 3$, and $\lambda > 0$. \par

(1) Let $u_+ \in H^1$. Then the NLS equation (1.1) (3.1) has a unique solution $u \in X_{loc}^1({I\hskip-1truemm R}) \cap X^1 ({I\hskip-1truemm R}^+)$ such that 
$$\parallel U(-t) \ u(t) - u_+ ; H^1 \parallel \to 0 \eqno(3.41)$$

\noindent when $t \to \infty$. \par

(2) Let $u \in X^1({I\hskip-1truemm R}^+)$ be a solution. Then there exists $u_+ \in H^1$ such that (3.41) holds.}\\

\noindent {\bf Sketch of proof.} The proof uses mainly Strichartz inequalities. In order to prove Part (1), one starts from the integral equation (3.3) with $u_+ = U(-t_0)u_0$ and $t_0 \to \infty$, namely 
$$u(t) = U(t)\ u_+ + i \int_t^{\infty} dt'\ U(t-t') \ g(\rho (t')) \ u(t') \eqno(3.42)$$

\noindent  and one solves that equation locally in a neighborhood of infinity in time, namely in $I = [T, \infty )$ for $T$ sufficiently large. The proof uses the estimate (3.4) followed by
$$\parallel u ; L^k(I, L^\ell )\parallel\ \leq \ C\parallel u ; L^{q_1}(I, \dot{H}_{r_1}^{\sigma_c})\parallel \ \leq \ C\parallel u ; X^1(I)\parallel \eqno(3.43)$$

\noindent for admissible $(q_1, r_1)$ with $k = q_1 < \infty$, namely
$$0 < 2/k = (n/2 - \sigma_c ) (1 - \delta ) = 2/q_1 = \delta (r_1)$$

\noindent which can always be realized for suitable $\delta$.\par

In order to prove Part (2), one estimates 
$$\parallel U(-t_1) \ u(t_1) - U(-t_2)\ u(t_2); H^1 \parallel\ = \ \parallel \int_{t_1}^{t_2} dt'\ U(t-t')\ g(\rho (t')) u(t'); H^1 \parallel$$
$$\leq \ C  \parallel u ; X^1(I)\parallel \  \parallel u ; L^{q_1}(I, \dot{H}_{r_1}^{\sigma_c})\parallel^{p-1} \eqno(3.44)$$

\noindent with $I = [t_1, t_2]$, and the last norm tends to zero when $t_1, t_2 \to \infty$ for $u \in X^1({I\hskip-1truemm R}^+)$ and $q_1 < \infty$. \par \hfill \sq \par

The previous proposition yields the existence of the wave operators in Part (1) and the fact that AC holds for solutions in $X^1({I\hskip-1truemm R})$ in Part (2). Putting together Propositions 3.1 and 3.2 yields AC for finite energy solutions.\\

\noindent {\bf Proposition 3.3.} {\it Let $0 < \sigma_c < 1$, $\sigma_c < 1/2$ for $n=1$, or equivalently $p > 1 + 4/n$, $p < 1 + 4/(n-2)$ for $n \geq 3$, and let $\lambda > 0$. Let $u$ be a finite energy solution of the NLS equation (1.1) (3.1), namely a solution $u \in X_{loc}^1({I\hskip-1truemm R})$. Then $u \in X^1({I\hskip-1truemm R})$ and there exist $u_\pm \in H^1$ such that 
$$\parallel U(-t) \ u(t) - u_\pm  ; H^1 \parallel \to 0 \eqno(3.45)$$
\noindent when $t \to \pm \infty$.} \\

We now turn to the Hartree equation (1.1) with $g$ given by (1.3). We assume that $V \in L^p$ for suitable $p$, for which we shall use the parameter $\sigma_c$ defined by
$$\sigma_c = n/2p-1 \ . \eqno(3.46)$$

\noindent The treatment extends in a trivial way to more general $V$ such as those considered in Proposition 2.2.\par

The main technical result is the following proposition.\\

\noindent {\bf Proposition 3.4.} {\it Let $n \geq 3$. Let $0 < \sigma_c < 1$, $\sigma_c \leq 1/2$ for $n = 3$, or equivalently $n/4 < p < n/2$, $p \geq 1$ for $n=3$. Let $V \in L^p$ be real even and let $u \in X_{loc}^1({I\hskip-1truemm R})  \cap L^{\infty} ({I\hskip-1truemm R}, H^1)$ be a solution of the Hartree equation (1.1) (1.3) such that $\rho = |u|^2 \in L^2({I\hskip-1truemm R} , \dot{H}^{(3-n)/2})$. Then $u \in X^1({I\hskip-1truemm R})$}. \\

\noindent {\bf Remark 3.5.} For nonnegative $V$, the Cauchy problem is globally well posed in $H^1$ and yields solutions $u \in X_{loc}^1({I\hskip-1truemm R})  \cap L^{\infty} ({I\hskip-1truemm R}, H^1)$.\\

\noindent {\bf Sketch of proof.} The proof follows the same pattern as that of Proposition 3.1. We start again from (3.3). Using the Strichartz estimates and the Young inequality, we estimate
$$\parallel u;X^1(I)\parallel\ \leq \ C \left ( \parallel u_0 ; H^1\parallel + \parallel V \parallel_p \ \parallel u; X^1(I)\parallel \ \parallel u; L^k(I, L^{\ell})\parallel^2 \right ) \eqno(3.47)$$

\noindent where now
$$\left \{ \begin{array}{l} 2/k = 1 - \delta \\ \\ n/\ell = n/2 - \sigma_c + \delta - 1\end{array} \right .$$

\noindent for some $\delta$, $0 \leq \delta \leq 1$. It is then sufficient to estimate $u \in L^k(I, L^{\ell})$ for $0 < 2/k \leq 1$ and
$$2/k + n/\ell = n/2 - \sigma_c\ .$$

\noindent The proof then proceeds as for the NLS equation. In particular one uses the estimate (3.8) with $k$, $\ell$ satisfying (3.9) so that $0 < \theta = 2(1- \delta ) \leq 1$, which ensures the condition $0 < 2/k \leq 1$. \par \hfill $\sq$ \par

The analogue of Proposition 3.2 can be proved for the Hartree equation \cite{11r} and the final result follows therefrom and from Propositions 3.4 and 2.2.\\

\noindent {\bf Proposition 3.5.}  {\it Let $n \geq 3$. Let $0 < \sigma_c < 1$, $\sigma_c \leq 1/2$ for $n = 3$, or equivalently $n/4 < p < n/2$, $p \geq 1$ for $n=3$. Let $V \in L^p$ be real radial and non increasing (and therefore nonnegative). Let $u$ be a finite energy solution of the Hartree equation (1.1) (1.3), namely a solution $u \in X_{loc}^1({I\hskip-1truemm R})$. Then $u \in X^1 ({I\hskip-1truemm R})$ and there exist $u_{\pm} \in H^1$ satisfying (3.45) when $t \to \pm \infty$.}

\section*{Appendix 1. Estimate of the RHS of (2.21) in H$^{\bf 1/2}$ for h = $|$x$|$} $\ \ \ $
\par\vskip 5 truemm

\noindent  {\bf Lemma} 
$$|<\rho, \nabla |x| \star {\rm Im} \ \overline{u} \nabla u > | \leq \ C \parallel u \parallel_2^2 \ \parallel  u; \dot{H}^{1/2}\parallel^2 \ . \eqno({\rm A}1.1)$$
\vskip 3 truemm
\noindent {\bf Proof.} We estimate 
$$|<\rho , \nabla |x| \star {\rm Im}\ \overline{u} \nabla u >|\ \leq \ |<\nabla u, u (\nabla |x| \star \rho)>|$$
$$\leq \ \parallel \omega^{1/2} u \parallel_2\ \parallel \omega^{1/2}u(\nabla |x| \star \rho ) \parallel_2$$
$$ \leq \ C \Big ( \parallel \omega^{1/2} u \parallel_2^2 \ \parallel \nabla |x| \star \rho \parallel_{\infty}$$
$$+ \  \parallel \omega^{1/2} u \parallel_2\ \parallel u \parallel_r \ \parallel  \omega^{1/2} (\nabla |x| \star \rho ) \parallel _{n/\delta } \eqno({\rm A}1.2)$$

\noindent with $\delta = \delta (r) > 0$, by fractional Leibniz inequalities. Clearly 
$$\parallel  \nabla |x| \star \rho  \parallel _{\infty} \ \leq \  \parallel \rho \parallel_1 \ . \eqno({\rm A}1.3)$$

\noindent We then use the fact that 
$$F(\nabla |x|) = C\ P \ \xi |\xi|^{-(n+1)}$$

\noindent where $F$ is the Fourier transform and $P$ denotes the principal value (\cite{18r}, Theorem~5, p.~73 with $k=1$) so that 
$$\omega^{1/2} (\nabla |x| \star \rho )  = C\ F^{-1} (\xi |\xi|^{-(n+1/2)} \widehat{\rho} (\xi )) = C\ x|x|^{-3/2} \star \rho$$

\noindent and therefore 
$$\parallel \omega^{1/2} (\nabla |x| \star \rho ) \parallel_{n/\delta} \ \leq \ C \parallel \rho \parallel_{s/2}\ = \ C \parallel u \parallel_s^2 \eqno({\rm A}1.4)$$ 

\noindent by the Hardy Littlewood Sobolev inequality (\cite{18r}, Theorem~1, p.~119), where
$$\delta + 2 \delta (s) = 1/2\ ,$$

\noindent provided $0 < \delta < n$ and $0 < \delta (s) < n/2$. The last term in (A1.2) is then estimated by
$$C  \parallel \omega^{1/2} u \parallel_2\ \parallel u \parallel_r \ \parallel u \parallel_s^2 \ \leq\ C \parallel \omega^{1/2} u \parallel_2^2 \ \parallel u \parallel_2^2  \eqno({\rm A}1.5)$$ 

\noindent by Sobolev inequalities, which together with (A1.2) and (A1.3) yields (A1.1). One can easily choose $r$ satisfying the required restrictions, for instance by taking $r=s$, which yields $\delta = 1/6$. \par \hfill $\sq$ \par

 \section*{Appendix 2. A quadratic identity for the NLKG equation}
 As mentioned in Section 2, the algebraic structure (2.1) (2.2) is realized for any system with symmetric conserved energy momentum tensor $T_{\lambda \mu}$, namely
$$\left \{ \begin{array}{l} \partial_t \sigma + \nabla \cdot j = 0 \\ \\ \partial_t j + \nabla \cdot T = 0 \end{array} \right . \eqno({\rm A}2.1)$$

\noindent where $\sigma = T_{00}$ is the energy density, $j_k = - T_{0k} = - T_{k0}$ is both the energy current and the momentum density, and $T = \{ T_{k \ell}\}$ is the space-space part of $T_{\lambda \mu}$. Here we use $\sigma$ instead of $\rho$ in order to keep the notation $\rho$ for $|u|^2$, greek (resp. latin) indices run from $0$ to $n$ (resp. from 1 to $n$), and the index $0$ refers to time. We shall also need the Minkowski metric $\eta_{\lambda \mu}$ with $\eta_{00} = - \eta_{jj} = 1$. We now consider the NLKG (or nonlinear wave NLW) equation
$$\sq u + g(\rho ) u = 0 \eqno({\rm A}2.2)$$

\noindent where $g = g(\rho )$ is a real function of $\rho = |u|^2$. The Lagrangian density is 
$$L(u) = |\partial_t u |^2 - |\nabla u|^2 - G(\rho ) \eqno({\rm A}2.3)$$

\noindent with $G$ defined by (2.9). The energy momentum tensor is well known to be
$$T_{\lambda \mu} = 2{\rm Re}\ \partial_{\lambda} \overline{u} \ \partial_{\mu} u  - \eta_{\lambda \mu} L \ . \eqno({\rm A}2.4)$$

\noindent We define as before for real even $h$
$$J = (1/2) \ < \sigma , h\star \sigma > \eqno({\rm A}2.5)$$

\noindent so that
$$M \equiv \partial_t  J = - \ <\sigma, \nabla h \star j>\eqno({\rm A}2.6)$$
$$\partial_t M = \partial_t^2 J = - \ <j, \nabla^2h \star j>\ +\ < \sigma, \nabla^2 h \star T>\ . \eqno({\rm A}2.7)$$

\noindent Substituting $\sigma = T_{00}$, $j_k = - T_{0k}$ and $T = \{ T_{k \ell}\}$ given by (A2.4) into (A2.7) yields
$$\partial_t M = - \ <2{\rm Re}\  \nabla \overline{u} \ \partial_t u, \nabla^2 h \star  2{\rm Re}\  \nabla \overline{u} \ \partial_t u>$$
$$+ \ < |\partial_t u|^2 + |\nabla u |^2 + G, \Delta h \star \left ( |\partial_t u|^2 - |\nabla u|^2 - G\right ) + \nabla^2 h \star 2 \nabla\overline{u} \nabla u>$$

\noindent and finally, ordering the terms by the powers of $G$
$$\partial_t M = - \ < G, \Delta h \star G>\ +\ <G, - 2 \Delta h \star |\nabla u|^2 + \nabla^2h \star 2\nabla\overline{u} \nabla u>$$
$$+ \ < |\partial_t u|^2 + |\nabla u |^2, \Delta h \star\left ( |\partial_t u|^2 - |\nabla u|^2\right ) + \nabla^2 h \star 2 \nabla\overline{u} \nabla u>$$
$$- \ <2{\rm Re}\  \nabla \overline{u}\  \partial_t u, \nabla^2 h \star  2{\rm Re}\  \nabla \overline{u} \ \partial_t u> \ . \eqno({\rm A}2.8)$$

\noindent In space dimension $n=1$, the linear term in $G$ vanishes and (A2.8) reduces to 
$$\partial_t M = - \ < G, h'' \star G>\ +\ < |\partial_t u|^2 + |\partial_x u|^2,  h'' \star (|\partial_t u|^2 + |\partial_x u|^2>$$
$$- \ <2\ {\rm Re}\  \partial_x \overline{u}\  \partial_t u, h'' \star  2\ {\rm Re}\  \partial_x  \overline{u}\  \partial_t u>$$
$$= \ < |\partial_t u - \partial_x u|^2, h'' \star |\partial_t u+ \partial_x u|^2>\  -\ < G,  h'' \star G> \eqno({\rm A}2.9)$$

\noindent so that if $h''$ has a given sign, the kinetic and potential terms have opposite signs, which precludes a straightforward use of that identity to derive estimates.

\section*{Appendix 3. Relation between the original and the quadratic Morawetz identities}
Here we rewrite the original Morawetz identity for the NLS equation in a form which exhibits its relation to the quadratic identity derived in Section~2. The original version starts from the quantity $<\nabla h, j>$ with $j = {\rm Im} \ \overline{u} \nabla u$. Using instead of $\nabla h$ a space translate thereof is then equivalent to consider the quantity
$$M_0 (x) = - \nabla h \star j \eqno({\rm A}3.1)$$

\noindent so that $M$ defined in (2.4) is simply $M = <\rho , M_0>$. Using the evolution equation of $j$ given by (2.2) (2.15) yields\par
\hskip 1.7 truecm $\partial_t M_0 = - \nabla^2 h \star T$
$$= \nabla^2 h \star ( \nabla \overline{u} \nabla u) + \Delta h \star (-(1/4) \Delta \rho + \rho g - G) \ . \eqno({\rm A}3.2)$$ 

\noindent For $h(x) = |x|$ and by the use of (2.22), we obtain
$$\partial_t M_0 = \nabla^2 h \star ( \nabla \overline{u} \nabla u) - (1/4) \Delta^2 |x| \star \rho + (n-1) |x|^{-1} \star (\rho g - G) \eqno({\rm A}3.3)$$ 

\noindent for $n > 1$, which is the original Morawetz identity for NLS. The first term in the RHS is nonnegative. The second term is positive for $n \geq 3$ only since
$$\Delta^2 |x| = \left \{\begin{array}{ll} - 8 \pi \delta (x) &\hbox{for $n=3$}\\ \\ - (n-1) (n-3)|x|^{-3} &\hbox{for $n\geq 4$}\end{array}\right .$$

\noindent Early applications of the method \cite{10r} \cite{12r} used (A3.3) to derive an estimate of the last term and were therefore restricted to space dimension $n \geq 3$. Taking the scalar product of $M_0$ with $\rho$ gives a useful quantity because $\rho$ satisfies the conservation law (2.1).

\section*{Appendix 4. Proof of (2.21) by regularization}
Let $u \in X_{loc}^1(I)$ be a solution of the NLS equation (2.7) and let $\varphi \in \mathscr{C}_0^{\infty} ({I\hskip-1truemm R}^n, {I\hskip-1truemm R}^+)$ with $\parallel \varphi \parallel _1 = 1$ denote a regularizing function of the space variable which will eventually converge to the measure $\delta$. We define $u_{\varphi} = \varphi \star u$, $\rho_\varphi = |u_\varphi |^2$ $j_\varphi = {\rm Im} \ \overline{u}_\varphi \nabla u_\varphi$, $f(u) = g(|u|^2)u$, $f_\varphi = \varphi \star f(u)$ and $f_{\not=} = f_\varphi - f(u_\varphi )$. Then $u_\varphi \in  \displaystyle{\mathrel{\mathop \cap_{\sigma \geq 0}}}\ \mathscr{C}(I, H^{\sigma})$ and $u_\varphi$ satisfies the equation
$$i \partial_t u_\varphi = - (1/2) \Delta u_\varphi + f_\varphi  \eqno({\rm A}4.1)$$ 

\noindent which is the regularized form of (2.7). From (A4.1), we obtain the regularized form of (2.1) (2.2) by direct computation, namely
$$\partial_t \ \rho_\varphi  + \nabla \cdot j_\varphi  = P_1 \ , \eqno({\rm A}4.2)$$
$$\partial_t \ j_\varphi + \nabla \cdot T(u_\varphi ) = P_2 \eqno({\rm A}4.3)$$

\noindent where
$$P_1 = 2\ {\rm Im} \left ( \overline{u}_\varphi f_{\not=}\right ) \ , \eqno({\rm A}4.4)$$ 
$$T_{k\ell} (u_{\varphi}) = {\rm Re} \nabla_k  \overline{u}_\varphi \nabla_{\ell} u_{\varphi} - \delta_{k\ell} \left ( (1/4) \Delta \rho_\varphi  - \rho_\varphi  g(\rho_\varphi ) + G(\rho_\varphi )\right ) \ ,  \eqno({\rm A}4.5)$$
$$P_2 = {\rm Re} \left (  f_{\not=} \nabla  \overline{u}_\varphi -  \overline{u}_\varphi \nabla f_{\not=} \right ) \ . \eqno({\rm A}4.6)$$

\noindent Using (A4.2) (A4.3) and the regularity properties of $u_\varphi$, we can derive the regularized version of (2.21) in the same way as in Section 2 (see (2.16)-(2.20)), namely
$$\int_{t_1}^{t_2} dt \Big \{ (1/2) \ <\nabla \rho_{\varphi}, \Delta h \star \nabla \rho_{\varphi}>\ + \ < \rho_{\varphi}, \Delta h \star \left (\rho_{\varphi} g(\rho_{\varphi}) - G(\rho_{\varphi})\right ) > $$
$$\left . + R(u_{\varphi} ) - S_{\varphi}\Big \} = - \ < \rho_{\varphi}, \nabla h \star j_{\varphi}> \right |_{t_1}^{t_2}\eqno({\rm A}4.7)$$

\noindent where
$$R(v) = \ < \overline{v}v, \nabla^2h \star \nabla \overline{v}\nabla v> \ - \ < \overline{v}\nabla v, \nabla^2 h \star \overline{v}\nabla v > \ , \eqno({\rm A}4.8)$$ 
\hskip 1.7 truecm $S_{\varphi} = -\ <j_{\varphi}, \nabla h \star P_1 >\ +\ < \rho_{\varphi} , \nabla h \star P_2>$
$$= 2 \ {\rm Re} \left \{ <u_\varphi \nabla \overline{u}_{\varphi} , \nabla h \star \overline{u}_\varphi f_{\not=}>\ + \ < \rho_{\varphi} , \nabla h \star  f_{\not=} \nabla \overline{u}_\varphi> \right \} \eqno({\rm A}4.9)$$

\noindent by a straightforward rewriting. \par

We now remove the cutoff, and without loss of generality, we restrict our attention to the case of a single power interaction of the form (3.1). We restrict ourselves to proving that
$$\lim_{\varphi \to \delta} \int_{t_1}^{t_2} dt\ S_{\varphi} = 0  \eqno({\rm A}4.10)$$

\noindent since the remaining terms in (A4.7) converge to the corresponding terms without $\varphi$ by estimates similar to (2.29) (2.31). We estimate
$$|S_{\varphi}| \ \leq \  4 \parallel u \parallel_2^2\  \parallel \nabla u \parallel_2 \ \parallel f_{\not=} \parallel_2 \eqno({\rm A}4.11)$$

\noindent since $\parallel \nabla h \parallel_\infty = 1$ by (2.22). From the identity
$$f_{\not=} = \varphi \star f(u) - f(u) + f(u) - f(u_{\varphi})$$

\noindent and from the estimates
$$|f(u)| \leq C|u|^p\ ,  \ |f(u) - f(u_{\varphi})| \leq C |u- u_{\varphi}| \left ( |u|^{p-1} + |u_{\varphi}|^{p-1} \right ) \ , $$

\noindent we obtain
$$\parallel  f_{\not=}  \parallel _2\ \leq \  \parallel \varphi \star f(u) - f(u) \parallel _2\ + \ C \parallel u \parallel _{2p}^{p-1} \  \parallel  u -  u_{\varphi}\parallel _{2p}\ . \eqno({\rm A}4.12)$$

\noindent We next estimate $\parallel u \parallel_{2p}$. For $n=1,2$ and for $n \geq 3$ with $p-1 \leq 2/(n-2)$, we estimate
$$\parallel u \parallel_{2p} \ \leq\ C \parallel u; H^1\parallel\ \in L_{loc}^{\infty}(I)\ .$$

\noindent For $n \geq 3$ and
$$2/(n-2) < p-1 < 4/(n-2)$$

\noindent we estimate
$$\parallel u \parallel_{2p} \ \leq \ C \parallel u; \dot{H}_r^1\parallel  \ \in L_{loc}^q(I)$$

\noindent where
$$0 < 2/q = \delta (r) = \delta (2p) - 1 < (n-2)/(n+2)$$

\noindent so that $\parallel f(u) \parallel_2 \ = C \parallel u \parallel_{2p}^p \ \in L_{loc}^k(I)$ with
$$0 < 2/k = (p-1)(n/2-1) - 1 < 1\ .$$

\noindent In both cases $\parallel  f_{\not=}  \parallel _2$ tends to zero for each $t \in I$ when $\varphi$ tends to $\delta$, and (A4.10) follows from (A4.11) and from the Lebesgue dominated convergence theorem.

\end{document}